\documentclass[12pt]{amsart}

\usepackage{amsfonts,latexsym,epsfig}

\def\R{\mathbb{R}}
\def\Z{\mathbb{Z}}
\def\d{\partial}

\theoremstyle{plain}

\begin{document}
\pagestyle{plain}

\title{Legendrian Mirrors and Legendrian Isotopy}
\author{Lenhard L.\ Ng}
\address{\newline \indent Department of Mathematics \newline
\indent Massachusetts Institute of Technology \newline
\indent 77 Massachusetts Avenue \newline
\indent Cambridge, MA 02139 \newline}
\email{lenny@math.mit.edu \newline 
\indent {\it URL}{\rm :} 
http://www-math.mit.edu/\~{}lenny/}
\renewcommand{\subjclassname}{\textup{2000} Mathematics Subject Classification}
\subjclass{Primary 53D12; Secondary 57M27, 57R17.}

\begin{abstract}
We resolve a question of Fuchs and Tabachnikov by showing that there
is a Legendrian knot in standard contact $\R^3$ with zero Maslov number
which is not Legendrian
isotopic to its mirror.  The proof uses the differential
graded algebras of Chekanov.
\end{abstract}

\date{First version: 20 April 2000; this version: 28 August 2000.}

\maketitle

A {\it Legendrian knot} in standard contact $\R^3$ is a knot which is 
everywhere tangent to the two-plane distribution induced by the contact
one-form $dz - y\,dx$.  Two Legendrian knots are {\it Legendrian isotopic}
if there is a smooth isotopy between them through Legendrian knots.
Fuchs and Tabachnikov \cite{FT} define an involution of Legendrian
knots as follows: given a Legendrian knot, let its {\it Legendrian mirror}
be the image of the knot under the diffeomorphism $(x,y,z) \mapsto
(x,-y,-z)$.  (This terminology is due to \cite{Che}; note, however,
that our standard contact form has the opposite orientation.)

Clearly a Legendrian knot and its mirror are isotopic as smooth knots; however,
this isotopy is not through contactomorphisms.  Thus there is no reason
{\it a priori} why a knot should be Legendrian isotopic to its mirror.
Nevertheless, an analysis of small examples shows that many knots are;
Fuchs and Tabachnikov (\cite{FT}, repeated in \cite{Tab} and \cite{Che})
ask whether this is true in general.  (They restrict
this question to knots with zero Maslov number, since mirroring negates
Maslov number; however, reversing knot orientation also negates Maslov
number, and so one could ask more generally if an oriented knot is always 
Legendrian isotopic to its mirror, with the opposite orientation if
necessary.)

We will show that the answer is negative by displaying a counterexample.
The proof uses the powerful differential graded algebra invariant
defined by Chekanov \cite{Che}.

\begin{figure}
\begin{center}
\epsfig{file=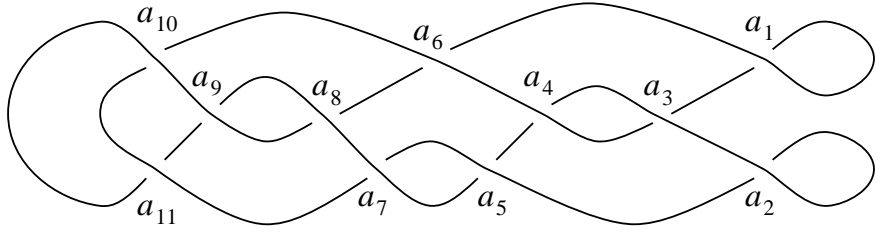,width=1.2in,angle=270}
\end{center}
\caption{Projection to the $xy$ plane of the Legendrian knot $K$, with 
crossings labelled.}
\label{sixtwo}
\end{figure}

Let $K$ be the (unoriented) Legendrian knot whose projection to the
$xy$ plane is given in Figure~\ref{sixtwo}.  (The diagram represents
a Legendrian knot because, e.g., it is the morsification of a
front diagram; see \cite{Fer}.)  This knot has smooth isotopy type
$6_2$, Maslov index 0, and Thurston-Bennequin invariant $-7$.
The Legendrian mirror $M(K)$ of $K$
is simply the reflection of this diagram about the $x$ (horizontal)
axis.  


\vspace{12pt}

\noindent 
\textbf{Proposition} \textit{The Legendrian mirrors $K$ and $M(K)$ are
not Legendrian isotopic.}

\begin{proof}
From \cite{Che}, it suffices to show that the Chekanov DGAs corresponding 
to $K$ and $M(K)$ have nonisomorphic graded homology algebras.
The Chekanov DGA associated to $K$ is the free associative unital
algebra $A(K)$ over $\Z/2$ generated by the crossings of $K$, which have
been labelled in Figure~\ref{sixtwo} as $a_1,\ldots,a_{11}$.
Of these generators, we calculate that $a_1,a_2,a_7,a_9,a_{10}$ have degree 1;
$a_3,a_4$ have degree 0; and $a_5,a_6,a_8,a_{11}$ have degree $-1$.
We can easily compute the differential on $A(K)$:
\[
\begin{array}{rclcrcl}
\d(a_1) &=& 1 + a_{10} a_5 a_3 & \hspace{0.25in} &
\d(a_6) &=& a_{11} a_8 \\
\d(a_2) &=& 1 + a_3 + a_3 a_6 a_{10} + a_3 a_{11} a_7 &&
\d(a_7) &=& a_8 a_{10} \\
\d(a_4) &=& a_5 + a_{11} + a_{11} a_7 a_5 + a_6 a_{10} a_5 &&
\d(a_9) &=& 1 + a_{10} a_{11} \\
\multicolumn{7}{c}{
\d(a_3) = \d(a_5) = \d(a_8) = \d(a_{10}) = \d(a_{11}) = 0.}
\end{array}
\]
The Chekanov DGA $A(M(K))$ associated to $M(K)$ has the same generators
and grading; its
differential is the same as the above differential, except with each
monomial reversed.

Note that 
$1 \in H_1(K) \cdot H_{-1}(K)$; that is, there exist $x \in H_1(K)$
and $y \in H_{-1}(K)$ such that $xy = 1 \in H_0(K)$.  (Choose $x=a_{10}$
and either $y=a_5 a_3$ or $y=a_{11}$.)  Since $A(M(K))$ is simply $A(K)$ with
monomials reversed, we conclude that 
$1 \in H_{-1}(M(K)) \cdot H_1(M(K))$.

We claim that $1 \not\in H_{-1}(K) \cdot H_1(K)$, which implies that
$H_*(K)$ and $H_*(M(K))$ are not isomorphic as graded algebras.
Suppose to the contrary that there exist $x,y,z \in A(K)$ of degree 
$-1,1,1,$ respectively, so that $1+\d z = xy$.  
Replace $a_3$ by $a_3+1$ and $a_{11}$ by $a_{11} + a_5$, as generators
of $A(K)$; the differential on $A(K)$ becomes
\[
\begin{array}{rclcrcl}
\d(a_1) &=& 1 + a_{10} a_5 + a_{10} a_5 a_3 & \hspace{0.25in} &
\d(a_6) &=& a_{11} a_8 + a_5 a_8 \\
\d(a_2) &=& a_3 + (1+a_3)(a_6 a_{10} + a_{11} a_7 + a_5 a_7) &&
\d(a_7) &=& a_8 a_{10} \\
\d(a_4) &=& a_{11} + a_{11} a_7 a_5 + a_6 a_{10} a_5 + a_5 a_7 a_5&&
\d(a_9) &=& 1 + a_{10} a_5 + a_{10} a_{11} \\
\multicolumn{7}{c}{
\d(a_3) = \d(a_5) = \d(a_8) = \d(a_{10}) = \d(a_{11}) = 0.}
\end{array}
\]
Then $1+xy=\d z$ is in the two-sided
ideal of $A(K)$ generated by $\{\d(a_i)\,:\, 1\leq i\leq 11\}$.

Let $A'$ be the free associative unital algebra over $\Z/2$ generated
by two variables $\alpha,\beta$, with grading which assigns degree $-1$ to
$\alpha$ and 1 to $\beta$.  We have a grading-preserving algebra map
$\pi\,:\,A(K)\rightarrow A'$ which sends $a_5$ to 
$\alpha$, $a_{10}$ to $\beta$, and all other generators to 0.
Write $x'=\pi(x)$ and $y'=\pi(y)$.
Then $1+x'y'$ is in the two-sided ideal $I$ of $A'$ generated
by $\{\pi(\d(a_i))\}$; by inspection, we see that $I$ is generated by the
single expression $1 + \beta \alpha$.

We conclude that $x'y' = 1$ in the algebra $A'/I$.  But
$A'/I$ is generated as a vector space by monomials of the form
$\alpha^i \beta^j$, where $i,j \geq 0$.  Since $x'$ has degree $-1$,
we can write $x' = \sum_{i\geq 1} c_i \alpha^i \beta^{i-1}$ (for
$c_i\in \Z/2$), and thus any monomial
$\alpha^i \beta^j$ in the expansion of $x'y'$ will satisfy $i \geq 1$.  
In particular, $x'y'$ cannot equal 1 in $A'/I$.  
This contradiction proves that $1 \not\in H_{-1}(K) \cdot H_1(K)$,
as desired.
\end{proof}

\vspace{16pt}

\section*{Acknowledgments}

The author would like to thank Tom Mrowka and Isadore Singer for their
encouragement and support.

\end{document}